%%%%%%%%%%%%%%%%%%%%  Oct 18, 1999 %%%%%%%%%%%%%%%%%%%%%%%%%

% begin.tex

%=========================================================================%
% load begin.tex only once, but keep count to match \bye commands
%=========================================================================%

\ifx\begin\undefined\else\global\advance\srcdepth by
1\expandafter \fi

\def\begin{}
\newcount\srcdepth
\srcdepth=1

\outer\def\bye{\global\advance\srcdepth by -1
  \ifnum\srcdepth=0
    \def\endcmd{\vfill\eject\nopagenumbers\par\vfill\supereject\end}
  \else\def\endcmd{}\fi
  \endcmd
}

%=========================================================================%
% initialize TeX
%=========================================================================%

%\magnification=\magstephalf
\magnification=\magstep 1
\baselineskip=13pt
%baselineskip=12pt
\hsize = 5.5truein
\hoffset = 0.5truein
\vsize = 8.5truein
\voffset = 0.2truein
\emergencystretch = 0.05\hsize

\overfullrule=0pt

\newif\ifblackboardbold

% comment out the following line if AMS msbm fonts aren't available
\blackboardboldtrue

%=========================================================================%
% select fonts
%=========================================================================%

\font\sectionfont=cmbx12

% Establish AMS blackboard bold fonts without using amssym.def, amssym.tex

\newfam\bboldfam
\ifblackboardbold
\font\tenbbold=msbm10
\font\sevenbbold=msbm7
\font\fivebbold=msbm5
\textfont\bboldfam=\tenbbold
\scriptfont\bboldfam=\sevenbbold
\scriptscriptfont\bboldfam=\fivebbold
\def\bbold{\fam\bboldfam\tenbbold}
\else
\def\bbold{\bf}
\fi

%=========================================================================%
% font size-changing command ("A Beginner's Book of TeX" p35, p275)
%=========================================================================%

\font\Arm=cmr8
\font\Ai=cmmi8
\font\Asy=cmsy8
\font\Abf=cmbx8
\font\Brm=cmr6
\font\Bi=cmmi6
\font\Bsy=cmsy6
\font\Bbf=cmbx6
\font\Crm=cmr5
\font\Ci=cmmi5
\font\Csy=cmsy5
\font\Cbf=cmbx5

\ifblackboardbold
\font\Abbold=msbm10 at 8pt
\font\Bbbold=msbm7 at 6pt
\font\Cbbold=msbm5
\fi

\def\smallmath{%
\textfont0=\Arm \scriptfont0=\Brm \scriptscriptfont0=\Crm
\textfont1=\Ai \scriptfont1=\Bi \scriptscriptfont1=\Ci
\textfont2=\Asy \scriptfont2=\Bsy \scriptscriptfont2=\Csy
\textfont\bffam=\Abf \scriptfont\bffam=\Bbf \scriptscriptfont\bffam=\Cbf
\def\rm{\fam0\Arm}\def\mit{\fam1}\def\oldstyle{\fam1\Ai}%
\def\bf{\fam\bffam\Abf}%
\ifblackboardbold
\textfont\bboldfam=\Abbold
\scriptfont\bboldfam=\Bbbold
\scriptscriptfont\bboldfam=\Cbbold
\def\bbold{\fam\bboldfam\Abbold}%
\fi
}

%=========================================================================%
% single-pass symbolic theorem labeling
%=========================================================================%

% Because this is a single-pass mechanism with no .aux file, forward
% references need to be declared in advance:

%   \forward{thm:main}{Theorem}{1.1}

% This is also the mechanism for "timely" declaration of labels, which
% will usually be buried within the corresponding theorem macros.
% A warning is issued if a label redeclaration is inconsistent, allowing
% forward references to be manually fixed.

%   \ref{thm:main} produces "Theorem~1.1"
%   \refs{thm:main} produces "Theorems~1.1"
%   \refn{thm:main} produces "1.1"

% Some TeX adapted from "The Advanced TeXbook" by David Salomon, chapter 9.

% Implementers: The code for \forward is subtle. Its second argument must
% be provided literally, e.g. "Theorem" rather that "\capitalize{theorem}".
% Its third argument must either be literal or a macro that expands
% directly to a literal, e.g. "\edef\numtoks{\number\proccount}".
% This use of \edef cannot be replaced by \def, which defers expansion.
% Failure to follow these rules will cause spurious warnings that forward
% references are inconsistent, when they are in fact consistent after
% expansion. Note the "Towers of Palo Alto" recreational math problem
% involving the iterated use of \expandafter to expand the first argument
% to \forwardsub before calling it.

\newlinechar=`@
\def\forwardmsg#1#2#3{\immediate\write16{@*!*!*!* forward reference should
be: @\noexpand\forward{#1}{#2}{#3}@}}
\def\nodefmsg#1{\immediate\write16{@*!*!*!* #1 is an undefined reference@}}

\def\forwardsub#1#2{\def\newref{{#2}{#1}}}

\def\forward#1#2#3{%
\expandafter\expandafter\expandafter\forwardsub\expandafter{#3}{#2}
\expandafter\ifx\csname#1\endcsname\relax\else%
\expandafter\ifx\csname#1\endcsname\newref\else%
\forwardmsg{#1}{#2}{#3}\fi\fi%
\expandafter\let\csname#1\endcsname\newref}

\def\firstarg#1{\expandafter\argone #1}\def\argone#1#2{#1}
\def\secondarg#1{\expandafter\argtwo #1}\def\argtwo#1#2{#2}

\def\ref#1{\expandafter\ifx\csname#1\endcsname\relax
  {\nodefmsg{#1}\bf`#1'}\else
  \expandafter\firstarg\csname#1\endcsname
  ~\expandafter\secondarg\csname#1\endcsname\fi}

\def\refs#1{\expandafter\ifx\csname#1\endcsname\relax
  {\nodefmsg{#1}\bf`#1'}\else
  \expandafter\firstarg\csname #1\endcsname
  s~\expandafter\secondarg\csname#1\endcsname\fi}

\def\refn#1{\expandafter\ifx\csname#1\endcsname\relax
  {\nodefmsg{#1}\bf`#1'}\else
  \expandafter\secondarg\csname #1\endcsname\fi}

%=========================================================================%
% widow control
%=========================================================================%

% usage:
% \widow{.2} % start new page if <.2 page left

\def\widow#1{\vskip 0pt plus#1\vsize\goodbreak\vskip 0pt plus-#1\vsize}

%=========================================================================%
% sections and theorems
%=========================================================================%

% use \showlabels or \showlabelsabove to display section and theorem labels

\def\marginlabel#1{}

\def\showlabelsabove{
\font\labelfont=cmss10 at 6pt
\def\marginlabel##1{\rlap{\smash{\raise 10pt\hbox{\labelfont##1}}}}
}

\newcount\seccount
\newcount\proccount
\seccount=0
\proccount=0

\def\stdskip{\vskip 9pt plus3pt minus 3pt}
\def\stdbreak{\par\removelastskip\penalty-100\stdskip}

\def\proof{\stdbreak\noindent{\sl Proof. }}

\def\qed{\vrule height 1.2ex width .9ex depth .1ex}

\def\Box{
  \ifmmode\eqno\qed
  \else\ifvmode\removelastskip\line{\hfil\qed}
  \else\unskip\quad\hskip-\hsize
    \hbox{}\hskip\hsize minus 1em\qed\par
  \fi\stdbreak\fi}

\def\ifempty#1#2\endB{\ifx#1\endA}
\def\makeref#1#2#3{\ifempty#1\endA\endB\else\forward{#1}{#2}{#3}\fi}

\outer\def\section#1 #2\par{
  \removelastskip
  \global\advance\seccount by 1
  \global\proccount=0\relax
                \edef\numtoks{\number\seccount}
  \makeref{#1}{Section}{\numtoks}
  \widow{.05}
  \vskip 24pt plus 6pt minus 6 pt
  \message{#2}
  \leftline{\marginlabel{#1}\sectionfont\numtoks\quad #2}
  \nobreak\stdskip}

\def\proclamation#1#2{
  \outer\expandafter\def\csname#1\endcsname##1 ##2\par{
  \stdbreak
  \advance\proccount by 1
  \edef\numtoks{\number\seccount.\number\proccount}
  \makeref{##1}{#2}{\numtoks}
  \noindent{\marginlabel{##1}\bf #2 \numtoks\enspace}
  {\sl##2\par}
  \stdbreak}}

\def\othernumbered#1#2{
  \outer\expandafter\def\csname#1\endcsname##1{
  \stdbreak
  \advance\proccount by 1
  \edef\numtoks{\number\seccount.\number\proccount}
  \makeref{##1}{#2}{\numtoks}
  \noindent{\marginlabel{##1}\bf #2 \numtoks\enspace}}}

\proclamation{definition}{Definition}
\proclamation{lemma}{Lemma}
\proclamation{proposition}{Proposition}
\proclamation{theorem}{Theorem}
\proclamation{corollary}{Corollary}
\proclamation{conjecture}{Conjecture}

\othernumbered{example}{Example}
\othernumbered{remark}{Remark}
\othernumbered{construction}{Construction}

%=========================================================================%
% enable postscript illustrations using epsf.tex
%=========================================================================%

% Usage:
% \draw{70}{fig}{} % draw fig.eps at 70% scale
% \draw{999}{fig}{} % draw fig.eps scaled to width of page

% Optional third argument can be multiple calls to \figtext; see below.
% More generally, the third argument is read in vertical mode, with the
% reference point at the lower left corner of the eps picture, whose
% dimensions are contained in the dimen registers \drawx and \drawy.
% This enables using TeX to generate the text that goes with the picture.
% To request that the picture be widened to respect the added text, 
% examine and modify the dimen registers \ngap, \egap, \sgap, \wgap.
% This is done automatically by the \figtext macro.

% These macros rely on "epsf.tex" which is the lowest level interface
% available for including encapsulated Postscript files in TeX documents.
% Rather that manually reading the .eps file to compute the nominal size,
% the \epsfbox macro is called twice, and two of its internal registers
% are examined after the first call. A major change to epsf.tex (unlikely)
% will require changes here. 

\input epsf

\newcount\figcount
\figcount=0
\newbox\drawing
\newcount\drawbp
\newdimen\drawx
\newdimen\drawy
\newdimen\ngap
\newdimen\sgap
\newdimen\wgap
\newdimen\egap

\def\drawbox#1#2#3{\vbox{
  \setbox\drawing=\vbox{\offinterlineskip\epsfbox{#2.eps}\kern 0pt}
  \drawbp=\epsfurx
  \advance\drawbp by-\epsfllx\relax
  \multiply\drawbp by #1
  \divide\drawbp by 100
  \drawx=\drawbp truebp
  \ifdim\drawx>\hsize\drawx=\hsize\fi
  \epsfxsize=\drawx
  \setbox\drawing=\vbox{\offinterlineskip\epsfbox{#2.eps}\kern 0pt}
  \drawx=\wd\drawing
  \drawy=\ht\drawing
  \ngap=0pt \sgap=0pt \wgap=0pt \egap=0pt 
  \setbox0=\vbox{\offinterlineskip
    \box\drawing \ifgridlines\drawgrid\drawx\drawy\fi #3}
  \kern\ngap\hbox{\kern\wgap\box0\kern\egap}\kern\sgap}}

\def\draw#1#2#3{
  \setbox\drawing=\drawbox{#1}{#2}{#3}
  \advance\figcount by 1
  \goodbreak
  \midinsert
  \centerline{\ifgridlines\boxgrid\drawing\fi\box\drawing}
  \smallskip
  \vbox{\offinterlineskip
    \centerline{Figure~\number\figcount}
    \smash{\marginlabel{#2}}}
  \endinsert}

\def\nextfigtoks{%
  \advance\figcount by 1%
  \edef\numtoks{\number\figcount}%
  \advance\figcount by -1}

\newif\ifgridlines
\newbox\figtbox
\newbox\figgbox
\newdimen\figtx
\newdimen\figty

\newdimen\bwd
\bwd=2sp % 2sp (1/32768") is smallest visible width for Textures

\def\hline#1{\vbox{\smash{\hbox to #1{\leaders\hrule height \bwd\hfil}}}}

\def\vline#1{\hbox to 0pt{%
  \hss\vbox to #1{\leaders\vrule width \bwd\vfil}\hss}}

\def\clap#1{\hbox to 0pt{\hss#1\hss}}
\def\vclap#1{\vbox to 0pt{\offinterlineskip\vss#1\vss}}

\def\hstutter#1#2{\hbox{%
  \setbox0=\hbox{#1}%
  \hbox to #2\wd0{\leaders\box0\hfil}}}

\def\vstutter#1#2{\vbox{
  \setbox0=\vbox{\offinterlineskip #1}
  \dp0=0pt
  \vbox to #2\ht0{\leaders\box0\vfil}}}

\def\crosshairs#1#2{
  \dimen1=.002\drawx
  \dimen2=.002\drawy
  \ifdim\dimen1<\dimen2\dimen3\dimen1\else\dimen3\dimen2\fi
  \setbox1=\vclap{\vline{2\dimen3}}
  \setbox2=\clap{\hline{2\dimen3}}
  \setbox3=\hstutter{\kern\dimen1\box1}{4}
  \setbox4=\vstutter{\kern\dimen2\box2}{4}
  \setbox1=\vclap{\vline{4\dimen3}}
  \setbox2=\clap{\hline{4\dimen3}}
  \setbox5=\clap{\copy1\hstutter{\box3\kern\dimen1\box1}{6}}
  \setbox6=\vclap{\copy2\vstutter{\box4\kern\dimen2\box2}{6}}
  \setbox1=\vbox{\offinterlineskip\box5\box6}
  \smash{\vbox to #2{\hbox to #1{\hss\box1}\vss}}}

\def\boxgrid#1{\rlap{\vbox{\offinterlineskip
  \setbox0=\hline{\wd#1}
  \setbox1=\vline{\ht#1}
  \smash{\vbox to \ht#1{\offinterlineskip\copy0\vfil\box0}}
  \smash{\vbox{\hbox to \wd#1{\copy1\hfil\box1}}}}}}

\def\drawgrid#1#2{\vbox{\offinterlineskip
  \dimen0=\drawx
  \dimen1=\drawy
  \divide\dimen0 by 10
  \divide\dimen1 by 10
  \setbox0=\hline\drawx
  \setbox1=\vline\drawy
  \smash{\vbox{\offinterlineskip
    \copy0\vstutter{\kern\dimen1\box0}{10}}}
  \smash{\hbox{\copy1\hstutter{\kern\dimen0\box1}{10}}}}}

\def\figtext#1#2#3#4#5{
  \setbox\figtbox=\hbox{#5}
  \dp\figtbox=0pt
  \figtx=-#3\wd\figtbox \figty=-#4\ht\figtbox
  \advance\figtx by #1\drawx \advance\figty by #2\drawy
  \dimen0=\figtx \advance\dimen0 by\wd\figtbox \advance\dimen0 by-\drawx
  \ifdim\dimen0>\egap\global\egap=\dimen0\fi
  \dimen0=\figty \advance\dimen0 by\ht\figtbox \advance\dimen0 by-\drawy
  \ifdim\dimen0>\ngap\global\ngap=\dimen0\fi
  \dimen0=-\figtx
  \ifdim\dimen0>\wgap\global\wgap=\dimen0\fi
  \dimen0=-\figty
  \ifdim\dimen0>\sgap\global\sgap=\dimen0\fi
  \smash{\rlap{\vbox{\offinterlineskip
    \hbox{\hbox to \figtx{}\ifgridlines\boxgrid\figtbox\fi\box\figtbox}
    \vbox to \figty{}
    \ifgridlines\crosshairs{#1\drawx}{#2\drawy}\fi
    \kern 0pt}}}}

% macros to add space to text on specified sides

\def\hpad#1#2#3{\hbox{\kern #1\hbox{#3}\kern #2}}
\def\vpad#1#2#3{\setbox0=\hbox{#3}\dp0=0pt\vbox{\kern #1\box0\kern #2}}

% macro to give one text string the apparent height of another

% macro to center one text string over another

\def\stack#1#2#3{\vbox{\offinterlineskip
  \setbox2=\hbox{#2}
  \setbox3=\hbox{#3}
  \dimen0=\ifdim\wd2>\wd3\wd2\else\wd3\fi
  \hbox to \dimen0{\hss\box2\hss}
  \kern #1
  \hbox to \dimen0{\hss\box3\hss}}}

% macros to hide size of trailing exponents

\def\hexp#1{%
  \setbox0=\hbox{${}^{#1}$}%
  \hbox to .5\wd0{\box0\hss}}

%=========================================================================%
% macros for matrices and arrows
%=========================================================================%

% typical usage:
%   \rightarrowmat{2pt}{4pt}{d & bd \cr \!-c & 0 \cr 0 & -ac \cr}

\def\bmatrix#1#2{{\smallmath\left[\vcenter{\halign
  {&\kern#1\hfil$##\mathstrut$\kern#1\cr#2}}\right]}}

\def\rightarrowmat#1#2#3{
  \setbox1=\hbox{\kern#2$\bmatrix{#1}{#3}$\kern#2}
  \,\vbox{\offinterlineskip\hbox to\wd1{\hfil\copy1\hfil}
    \kern 3pt\hbox to\wd1{\rightarrowfill}}\,}

\def\leftarrowmat#1#2#3{
  \setbox1=\hbox{\kern#2$\bmatrix{#1}{#3}$\kern#2}
  \,\vbox{\offinterlineskip\hbox to\wd1{\hfil\copy1\hfil}
    \kern 3pt\hbox to\wd1{\leftarrowfill}}\,}

\def\rightarrowbox#1#2{
  \setbox1=\hbox{\kern#1\hbox{\smallmath #2}\kern#1}
  \,\vbox{\offinterlineskip\hbox to\wd1{\hfil\copy1\hfil}
    \kern 3pt\hbox to\wd1{\rightarrowfill}}\,}

\def\leftarrowbox#1#2{
  \setbox1=\hbox{\kern#1\hbox{\smallmath #2}\kern#1}
  \,\vbox{\offinterlineskip\hbox to\wd1{\hfil\copy1\hfil}
    \kern 3pt\hbox to\wd1{\leftarrowfill}}\,}

%=========================================================================%
% quire macros for preview mode and making booklets
%=========================================================================%

% \legalbooklet{20} makes a booklet from legal paper in landscape
% orientation, where "20" is the page count. To preview, give a negative
% pagecount. Either print using the legal duplex option on a modern laser
% printer, or struggle to simulate this effect manually. Bind using a long
% reach stapler.

% \preview squeezes two pages side by side in landscape orientation. It
% is not suitable for printing, but ideal for previewing on a two page
% monitor.

% \twoup squeezes two pages onto letter paper in landscape mode,
% suitable for printing.

% Each of these macros calls the file "quire.tex"

\def\bookletdims{
  \hsize=5.25truein
  \vsize=7truein
}

\def\legalbooklet#1{
  \input quire
  \bookletdims
  \htotal=7.0truein
  \vtotal=8.5truein
  % below computed from above
  \hoffset=\htotal
  \advance\hoffset by -\hsize
  \divide\hoffset by 2
  \voffset=\vtotal
  \advance\voffset by -\vsize
  \divide\voffset by 2
  \advance\voffset by -.0625truein
  \shhtotal=2\htotal
  % below doesn't need to change
  \horigin=0.0truein
  \vorigin=0.0truein
  \shstaplewidth=0.01pt
  \shstaplelength=0.66truein
  \shthickness=0pt
  \shoutline=0pt
  \shcrop=0pt
  \shvoffset=-1.0truein
  \ifnum#1>0\quire{#1}\else\qtwopages\fi
}

\def\preview{
  \input quire
  \bookletdims
  \hoffset=0.1truein
  \vtotal=8.5truein
  \shhtotal=14truein
  % below computed from above
  \voffset=\vtotal
  \advance\voffset by -\vsize
  \divide\voffset by 2
  \advance\voffset by -.0625truein
  \htotal=2\hoffset
  \advance\htotal by \hsize
  % below doesn't need to change
  \horigin=0.0truein
  \vorigin=0.0truein
  \shstaplewidth=0.5pt
  \shstaplelength=0.5\vtotal
  \shthickness=0pt
  \shoutline=0pt
  \shcrop=0pt
  \shvoffset=-1.0truein
  \qtwopages
}

\def\twoup{
  \input quire
  \hsize=4.79452truein % 5.25/1.095
  \vsize=7truein
  \vtotal=8.5truein
  \shhtotal=11truein
  % below computed from above
  \hoffset=-2\hsize
  \advance\hoffset by \shhtotal
  \divide\hoffset by 6
  \voffset=\vtotal
  \advance\voffset by -\vsize
  \divide\voffset by 2
  \advance\voffset by -12truept
  \htotal=2\hoffset
  \advance\htotal by \hsize
  % below doesn't need to change
  \horigin=0.0truein
  \vorigin=0.0truein
  \shstaplewidth=0.01pt
  \shstaplelength=0pt
  \shthickness=0pt
  \shoutline=0pt
  \shcrop=0pt
  \shvoffset=-1.0truein
  \qtwopages
}

%=========================================================================%
% timestamp (adapted from eplain.tex)
%=========================================================================%

\newcount\countA
\newcount\countB
\newcount\countC

\def\monthname{\begingroup
  \ifcase\number\month
    \or January\or February\or March\or April\or May\or June\or
    July\or August\or September\or October\or November\or December\fi
\endgroup}

\def\dayname{\begingroup
  \countA=\number\day
  \countB=\number\year
  \advance\countA by 0 % adjust after each leap day
  \advance\countA by \ifcase\month\or
    0\or 31\or 59\or 90\or 120\or 151\or
    181\or 212\or 243\or 273\or 304\or 334\fi
  \advance\countB by -1995
  \multiply\countB by 365
  \advance\countA by \countB
  \countB=\countA
  \divide\countB by 7
  \multiply\countB by 7
  \advance\countA by -\countB
  \advance\countA by 1
  \ifcase\countA\or Sunday\or Monday\or Tuesday\or Wednesday\or
    Thursday\or Friday\or Saturday\fi
\endgroup}

\def\timename{\begingroup
   \countA = \time
   \divide\countA by 60
   \countB = \countA
   \countC = \time
   \multiply\countA by 60
   \advance\countC by -\countA
   \ifnum\countC<10\toks1={0}\else\toks1={}\fi
   \ifnum\countB<12 \toks0={\sevenrm AM}
     \else\toks0={\sevenrm PM}\advance\countB by -12\fi
   \relax\ifnum\countB=0\countB=12\fi
   \hbox{\the\countB:\the\toks1 \the\countC \thinspace \the\toks0}
\endgroup}

\def\timestamp{\dayname, \the\day\ \monthname\ \the\year, \timename}

%==========================================================================
% macros (specific to this paper)
%==========================================================================

% surround with $ $ if not already in math mode
\def\enma#1{{\ifmmode#1\else$#1$\fi}}

%\showlabelsabove
\othernumbered{question}{Question}
\othernumbered{algorithm}{Algorithm}
\othernumbered{example}{Example}

\forward {modulebounds}{Section}{4}
\forward {cohomology bounds}{Section}{3}
\forward {support of coho}{Section}{1}
\forward {dependence for monomial}{Theorem}{1.1}
\forward {local cohomology}{Section}{2}
\forward {Huneke}{Example}{5.4}
\forward {alg}{Section}{5}
\forward {support of coho global}{Section}{2}
\forward {problem on lc}{Section}{6}
\input diagrams.tex
\input epsf.tex

\def \Z {{\bf Z}}
\def \iso {\cong}
\def \P {{\bf P}}

\def \O {{\cal O}}

\def \ZZ {{\bf Z}}
\def \RR {{\bf R}}
\def \Ch {{\bf D}}
\def \D  {{\bf D}}
\def \th {{\rm th}}

\def \O {{\cal O}}

\def \Tor {\mathop{\rm Tor}\nolimits}
\def \Ext {\mathop{\rm Ext}\nolimits}
\def \H {\mathop{\rm H}\nolimits}

\def \Hom {\mathop{\rm Hom}\nolimits}

\def \codim {\mathop{\rm codim}\nolimits}
\def \neg {\mathop{\rm neg}}
\def\fix#1{{\bf ((**** #1 ****))}}
\def \mustata {Musta\c t\v a}

\centerline{\bf Cohomology On Toric Varieties}
\centerline{\bf and}
\centerline{\bf Local Cohomology With Monomial Supports}
\bigskip
\centerline{David Eisenbud, Mircea \mustata, and Mike
Stillman\footnote{$^{*}$}{\rm The authors
are grateful to the NSF for partial
support during the preparation of this work.}}
\bigskip

In this note we describe aspects of the cohomology of coherent sheaves
on a complete toric variety $X$ over a field $k$ and, more generally,
the local cohomology, with supports in a monomial ideal, of a finitely
generated module
over a polynomial ring $S$. This leads to an efficient way of computing such
cohomology, for which we give explicit algorithms.

The problem is finiteness.
The $i^\th$ local cohomology of an $S$-module $P$
with supports in an ideal $B$ is the limit
$$
\H_B^i(P)=\lim_{\longrightarrow\atop \ell}\Ext^i(S/B_\ell, P),
$$
where $B_\ell$ is any sequence of ideals that is cofinal with the
powers of $B$. We will be interested in the case where $S$ is a
polynomial ring, $P$ is a finitely generated module, and $B$ is a
monomial ideal.  The module on the left of this equality is almost
never finitely generated (even when $P=S$), whereas the module
$\Ext^i(S/B_\ell, P)$ on the right is finitely generated, so that
the limit is really necessary.

We can sometimes restore finiteness by considering the homogeneous
components of $H^i_B(P)$ for a suitable grading. For example, in
the case where $B$ is a monomial ideal and $P=S$
the modules on
both sides are $\Z^n$ graded, and for each $p\in\Z^n$
the component $H^i_B(S)_p$ is a finite
dimensional vector space.  From this and the Hilbert Syzygy Theorem,
one sees easily that the corresponding finiteness
holds for any finitely generated
$\Z^n$-graded module $P$.

To get a computation that works for more than
$\Z^n$-graded modules, we work with coarser gradings;
that is we consider a grading in an abelian group $\D$ that
is a homomorphic image of $\Z^n$.
Of course if the grading is too coarse, we will lose finiteness
again. In this paper we give a criterion on the
grading for such finiteness to hold. When it
holds we study the convergence of the limit above,
and show how to compute,
for each $\delta\in\D$, an explicit
ideal $B^{[\ell]}$ such that the natural map
$$
\Ext^i(S/B^{[\ell]}, P)_\delta\rTo \H_B^i(P)_\delta
$$
is an isomorphism. The computation involves the solution of
a linear programming problem that depends on
the grading $\D$ and on the syzyies of
$P$.

The results in this paper were motivated by the problem
of computing the cohomology of coherent sheaves on a complete
toric variety, and such varieties provide the
most interesting cases for which our finiteness condition
is satisfied. The connection is via David Cox's [1995] notion of the
homogeneous coordinate ring of a toric variety,
described in \ref{support of coho global}.
The homogeneous coordinate ring is a polynomial ring
$S=k[x_1,\dots,x_n]$ equipped with a grading in an abelian group
$\D$ (the group of invariant divisor classes on $X$) and an
{\it irrelevant ideal} $B=B_X$ generated by square-free monomials,
defined from the fan associated to $X$. The data $B$ and $\D$
satisfy our finiteness condition.

Folowing ideas of Cox and \mustata \ explained in
\ref{support of coho global}
we may represent any coherent sheaf on a toric variety $X$
by giving a $\D$-graded module on a polynomial ring $S$,
and the cohomology of the sheaf is given by a formula
similar to that for local cohomology above. Thus we get an
explicit computation of sheaf cohomology in this case.

In the first section of this paper we treat various general
results on local cohomology with monomial supports. We show that
the local cohomology $\H^i_B(S)_p$, for $p\in \Z^n$,
depends only on which coordinates of $p$ are negative.
We analyze the condition that for a coarser grading $\Z^n\to \D$
the homogeneous components  $\H^i_B(S)_\delta$ are finite
dimensional for all $\delta\in\D$ --- this condition is
satisfied for example when $B$ and $\D$ come from the homogeneous
coordinate ring of a complete toric variety. In the second section
we translate these results to the case of toric varieties, and
also give a new topological charaterization of the $p\in\Z^n$ for
which  $\H^i_B(S)_p\neq 0$.

In section 3, assuming finiteness,
we determine $\ell$ such that the
natural map
$
\Ext^i(S/B^{[\ell]}, S)_\delta \rTo H_B^i(S)_\delta
$
is an isomorphism,
and similarly for sheaf cohomology.
In section 4 we use syzygies to extend this to all
modules $P$.

Section 5 is devoted to algorithms made from these results.

In section 6 we present some basic problems, partially solved in this
paper in the monomial case: when are the maps
$
\Ext^i(S/I, S)\rTo \H_I^i(S)
$
monomorphisms? In general, given an ideal $B$,
we do not even know that there are
ideals $I$ with the same radical as $B$ such that this
map is a monomorphism, and we
display a method proposed by Huneke that gives a
non-existence criterion.

Using the computation of sheaf cohomology via local cohomology and the
results above, \mustata [1999b] has proved a cohomology vanishing
result on toric varieties including one that strengthens a version of
the well-known vanishing theorem of Kawamata and Viehweg in this case,
and is valid in all characteristics. The main point is again the fact that
$H_B^i(S)_p$, for $p\in\Z^n$ depends only on which coordinates
of $p$ are negative. His result is:

\theorem{} Let $X$ be an arbitrary toric variety and $D$ an invariant Weil
divisor. If there is $E=\sum_{j=1}^da_jD_j$,with $a_j\in {\bf Q}$,
$0\leq a_j\leq 1$ and $D_j$ prime invariant Weil divisors such that for
some integers $m\geq 1$ and $i\geq 0$, $m(D+E)$ is integral and Cartier
and 
$$
\H^i(\O_X(D+m(D+E)))=0,
$$
then $\H^i(\O_X(D))=0$. In particular, if $X$ is complete and there is $E$
as before such that $D+E$ is {\bf Q}-ample, then $\H^i(\O_X(D))=0$,
for all $i\geq 1$.

An important step toward the theorems of this paper is a
result of
\mustata \  [1999a], motivated by this project. It shows that the local
cohomology of $M=S$ with monomial supports
can be computed as a union, not just
as a limit, of suitable $Ext$ modules.
For the reader's convenience, we state it here.
 We write $B^{[\ell]}$ for
the ideal generated by the $\ell^\th$ powers of monomials in $B$.

\theorem{lc as union} Let $B\subset S = k[x_1, \ldots, x_n]$
be a square-free monomial ideal.
For each integer $j\geq 0$ the natural map
$$
\Ext^j(S/B^{[\ell]}, S) \longrightarrow \H^j_B(S)
$$
is an injection, and its image is the submodule of
$\H^j_B(S)$ consisting of all elements of degree
$ \geq (-\ell,-\ell,\ldots,-\ell)$.\Box

It would be nice to have such results for more general ideals
$B$; see \ref{problem on lc} for some remarks on this problem.

We are grateful to Louis Billera, Craig Huneke,
Sorin Popescu and Greg Smith for helpful
discussions regarding parts of this material.

Henceforward in this paper $k$ denotes an arbitrary field,
and $B$ denotes a reduced monomial ideal in a polynomial
ring $S=k[x_1, \dots ,x_n]$.

\section {support of coho} The support of $\H^i_B(S)$

In this section we establish our basic results for the local
cohomology of the ring with supports in a monomial ideal $B$.
Since the local cohomology depends only on the radical of $B$,
we may assume that $B$ is generated by square-free monomials.
In the next section we explain the parallel theory
for the cohomology of sheaves on a toric variety, with some
refinements possible only in that case.

As $B$ is a monomial ideal, the local cohomology $\H^i_B(S)$
is naturally $\Z^n$-graded. Our main goal is to describe the set of
indices $p\in\Z^n$ for which $\H^i_B(S)_p\neq 0$. We first show that
this condition depends only on which components of $p$ are
negative. Set $\neg (p)=\{i\in\{1,\dots,n\}\mid p_i<0\}$.

\theorem{dependence for monomial}
If $p$, $q\in\Z^n$ satisfy $\neg(p)=\neg(q)$, then there is a canonical
isomorphism
$$
\H_B^i(S)_p\,\cong\,\H_B^i(S)_q.
$$

\proof
We use the computation of local cohomology as
\v Cech cohomology.
Given a sequence of elements ${\bf f} = \{f_1, \ldots, f_r\}$ of
a ring $R$, and an $R$-module $P$, the
{\it \v Cech complex\/} $C({\bf f},P)$
is the complex with
$$
C^i({\bf f},P) = \bigoplus_{j_1 < \ldots <
j_i} P_{f_{j_1} f_{j_2}\ldots f_{j_i}},
$$
and differential
$$
\partial^i : C^i({\bf f}, P) \longrightarrow C^{i+1}({\bf f},
P)
$$
which is the alternating sum of the localization maps. The
$i^\th$ cohomology group of this complex is the $i^\th$
local cohomology of $P$ with supports in the ideal generated
by the $f_i$; see for example Brodmann and Sharp [1998].

Let $(m_1, \ldots, m_r)$ be monomial generators for
the ideal $B$, and write {\bf m} for the sequence $m_1,\dots,m_r$.
Given $p \in \Z^n$, let $C({\bf f},S)_p$ denote the
complex of vector spaces that is the degree $p$ part of the
\v Cech complex
$C({\bf m},S)$. Let $m_J$ be the least common multiple of
the monomials $m_j$, $j\in J$. It is easy to check that
$$
C({\bf m},S)_p = \bigoplus_{\{J \mid \neg(p) \subset supp(m_J)\}} k,
$$
with differential mapping the $J^{th}$ component to the $J'^\th$
component equal to zero unless $J'$ has the form $J'=J\cup{j}$,
while in this case it is $(-1)^e$ where $e$ is the position of $j$
in the set $J'$.

Thus the complex $C({\bf m},S)_p$ depends only on $\neg(p)$.
The homology
of this complex is
$
\H_B^i(S)_p.
$
\Box

In view of \ref{dependence for monomial} it is useful to define
$$
L_I=\{ p\in\Z^n \mid \neg(p)=I \}.
$$
If we write $Supp(M)$ for the set of $\Z^n$-degrees
in which a $\Z^n$-graded module $M$ is nonzero,
then part of
\ref{dependence for monomial} asserts that $Supp(\H^i_B(S))$
is a union of certain $L_I$. We set
$$
\Sigma_i=\{ I \subset \{1,\dots,n \} \mid \H^i_B(S)_p\neq 0\,
{\rm for}\, p\in L_I \}.
$$

It will also be useful to define $C_I \subset \ZZ^n$ to be the orthant
$$
C_I =\bigl\{p\in\ZZ^n\mid  \cases{p_i \leq 0 &if $j \in I$\cr
                             p_i \geq 0 &if $j \notin I$}
     \bigr\}
$$
and
$$
(p_I)_j=\cases{-1 &if $j \in I$\cr
                0 &if $j \notin I$.}
$$
Note that $L_I=p_I+C_I\subset C_I$.

If $\Delta$ is a simplicial complex on $\{1,\ldots,n\}$, the
Alexander dual $\Delta^*$ is defined to be the simplicial complex
  $$\Delta^* := \{ \{1,\ldots,n\} \setminus f \mid f \not\in
\Delta \},$$
and the Alexander dual of a square-free monomial ideal $B$ is the
corresponding notion for Stanley-Reisner ideals:  $B^*$ is the
monomial ideal generated by the square-free monomials in
  $$(x_1^2, \ldots, x_n^2) : B.$$

Part (c) of the following corollary seems to be the most efficient way to
compute
the $\Sigma_i$.

\corollary{find Sigma} The following are equivalent.\hfil\break\indent
(a) $I\in\Sigma_i$ \hfil\break\indent
(b)$\Ext^i(S/B,S)_{p_I}\neq 0$ \hfil\break\indent
(c)$\Tor^S_j(S/B^*, k)_{-p_I} \neq 0$, where $j = \#I - i + 1$.

\proof
By \mustata's theorem \ref{lc as union} we have
$\Ext^i(S/B,S)_{p_I}=\H^i_B(S)_{p_I}$
and \ref{dependence for monomial} gives the equivalence of (a) and (b).
The equivalence of (b) and (c) is Corollary 3.1 in  \mustata[1999a] \Box

Now we introduce a grading coarser than the $\Z^n$-grading.
Let
$$
0\to M \rTo^\rho \Z^n \rTo^\phi \D\to 0
$$
be an exact sequence of abelian groups. Any
$\Z^n$-graded module $N$ can be regarded as a $\D$-graded
module by setting $N_\delta=\oplus_{p\in\phi^{-1}\delta}N_p$
for each $\delta\in\D$.
We are interested in the finiteness of
$\H_B^i(S)_\delta$ for all $\delta\in \D$.
From the \v Cech complex (or from \ref{lc as union}) we see
that $\H^i_B(S)_p$ is a finite dimensional vector space
for all $p\in\Z^n$. Thus we have

\corollary{empty intersection}
$(p+M)\cap L_I$ is finite
for every $I \in \Sigma_i$ and every $p\in\Z^n$, iff
$\H^i_B(S)_{\delta}$ is a finite dimensional vector space
for every $\delta\in\D$.\Box

The finiteness of the components $\H^i_B(S)_\delta$
translates into a condition on the subsets $I$ in the $\Sigma_i$
and the subgroup $M\subset \Z^n$ defining the coarse grading:

\proposition{inter equivs} With notation as above, let $I'$ be the
complement of $I$.
The following are equivalent:
\item {(a)} For every $p\in\Z^n$, the set $(p+M)\cap L_I$
is finite.
\item{(b)} $M\cap C_I=0$.
\item{(b')} $M\cap C_{I'}=0$
\item{(c)} For any $p, q, r\in \ZZ^n$
such that  $q\neq r\in(p+M)\cap C_I$,
there are indices $i,j$ such that $|q_i|>|r_i|$
and $|q_j|<|r_j|$.
\item{(d)} For every $p\in \ZZ^n$
 the set $(p+M)\cap C_I$ is finite.
\item{(d')} For some $p\in\Z^n$, the set $(p+M)\cap C_I$
is finite and nonempty.

\proof
$(a)\Rightarrow (b)$: Suppose on the contrary that
$0\neq p\in C_I\cap M$. Choose $q\in L_I$.
Then $q+s p\in (q+M)\cap L_I$ for every positive integer $s$,
showing that $(q+M)\cap L_I$ is not finite.
\hfill\break
$(b)\equiv (b')$: $M\cap C_I = -(M\cap C_{I'})$.
\hfill\break
$(b)\Rightarrow (c)$: If all the coordinates of $q$
had absolute value $\geq$ the corresponding coordinates
of $r$ then $q-r\in M\cap C_I$, and similarly
for the other inequality.
\hfill\break
$(c)\Rightarrow (d)$: Suppose contrary to $c)$ that
$(p+M)\cap C_I$ is infinite.
Any set of elements in an orthant contains a finite set
of minimal elements with respect to the partial order by absolute values of
the coordinates, and every element is $\geq$ to one of these in
the partial order.
(Proof: Moving to the first quadrant we may think of the
elements as monomials in the polynomial ring. By the Hilbert
Basis Theorem, a finite subset generates the same ideal
as the whole set.) Thus two elements of
$(p+M)\cap C_I$ are comparable, contradicting $b)$.
\hfill\break
$(d)\Rightarrow (d')$: It is enough to pick $p\in C_I$,
so that $(p+M)\cap C_I$ is nonempty.
\hfill\break
$(d')\Rightarrow (b)$: If $q\in(p+M)\cap C_I$ and
$r\in M\cap C_I$ then
$q+r\in(p+M)\cap C_I$.
Thus if $(p+M)\cap C_I$ is finite then $M\cap C_I$ is finite.
As the latter is a cone, it must be 0.
\hfill\break
$(d)\Rightarrow (a)$: This is obvious since $L_I\subset C_I$.
\Box

Of particular importance to us is the image
of the cones $L_I$ (with $I\in \Sigma_i$) in $\D$. The following
result is what ultimately allows us to understand the convergence
of  $\Ext^i(S/B^{[\ell]},S)_\delta$ to $\H^i_B(S)_\delta$ for
$\delta\in\D$.

\corollary{coarse image} If the condition in \ref{empty intersection}
is satisfied, then
for $I\in \Sigma_i$ the projection
$\phi(C_I)$ is a pointed cone in the sense that if $x,-x\in\phi(C_I)$ then
$x=0$. Further, $p_I$ maps to a nontorsion element of $D$.

\proof It is easy to see that
$\phi(C_I)$ contains both a nonzero vector and its negative
if and only if the kernel of $\phi$ meets an open face of
the cone $C_I$ without containing it. By part (b) in \ref{inter equivs}
the kernel $M$ meets only the zero face.
\Box

\section{support of coho global} The case of toric varieties

Let $X$ be a toric variety
which is always assumed to be nondegenerate,
meaning that the corresponding fan is not contained
in a hyperplane.  
 In this section we
will introduce the technique we will use to handle sheaves on $X$,
and we describe the
set of torus invariant divisors $D$ for which
$\H^i(\O_X(D))$ is nonzero..

We begin by reviewing Cox's homogeneous coordinate ring,
and the basic computation of cohomology that it allows.
Let
$\Delta$ be the fan in $\Z^d$ corresponding to $X$
(see for example Fulton [1993] for notation and terminology),
and suppose that $x_1,\dots,x_n$ are the edges (one-dimensional
cones) of $\Delta$. The torus-invariant divisor classes correspond
to the elements of the cokernel $\D$ of the map
represented by a matrix whose rows are the coordinates of the
first integral points of the $n$ edges of $\Delta$, so that we have
an exact sequence
$$
M:=\Z^d\rTo^\rho \Z^n\rTo^\phi \D.
$$
As before, the map $\phi$ defines a grading by $\D$ on the polynomial ring
$S := k[x_1,\dots,x_n]$

Cox [1995] defines the {\it homogeneous coordinate ring\/}
of $X$ to be
the polynomial ring $S$
together with the $\D$-grading and the {\it irrelevant ideal\/}
$$
B\ =\  (\{\prod_{x_i \not\in \sigma} x_i \mid \sigma \in \Delta\} ).
$$

Given a $\D$-graded $S$-module $P$, Cox constructs a quasicoherent
sheaf $\tilde P$ on $X$ by localizing just as in the case of
projective space. Coherent sheaves come from finitely generated
modules:

\theorem{what is a sheaf} Every coherent
$\O_X$ module may be written as $\tilde{P}$, for a finitely
generated $\D$-graded $S$-module $P$.

A proof is given by Cox [1995]  in the simplicial
case and in \mustata [1999b] in general. \Box

For any $\D$-graded $S$-module $P$ and any $\delta\in \D$
we may define $P(\delta)$ to be the graded module with components
$P(\delta)_\epsilon=P_{\delta+\epsilon}$ and we set
$$
\H^i_*(\tilde P)=\oplus_{\delta\in \D}\H^i(\widetilde {P(\delta)}).
$$

We have $\H^0(\O_X(\delta))=S_\delta$ for
each $\delta\in \D$. In general we write
$$
\H^i_*(\O_X)=\oplus_{\delta\in \D}\H^i(\O_X(\delta))
$$
so that $\H^0_*(\O_X) = S$. In fact each $\H^i_*(\O_X)$ is a $\Z^n$-graded
$S$-module.
We can compute
$H^i_*(\tilde{P})$ in terms of
local cohomology or limits of Ext modules.
The results and proofs are easy generalizations of results in
Grothendieck [1967].

We will consider also the shifted \v Cech complex $C_{[1]}({\bf f},P)$,
where $C_{[1]}^i({\bf f}, P)=C^{i+1}({\bf f}, P)$ for $i\geq 0$
and $C_{[1]}^{-1}({\bf f}, P)=0$.
If $P$ and the $f_i$ are all $\D$-graded
then each map in these complexes is
graded of degree zero, and so the cohomology modules will also be
$\D$-graded.

\theorem{}
Let $P$ be a $\D$-graded $S$-module. Let $\tilde{P}$ be the
corresponding quasi-coherent sheaf on $X$ as above. If
the irrelevant ideal $B\subset S$ is generated by
the monomials $m_i$ then
$$\H^i_*(\tilde{P}) \,\cong\, H^i(C_{[1]}({\bf m}, P)),$$
as graded $S$-modules.

The proof is immediate from the definition of $\tilde{P}$
and the computation of cohomology as \v Cech cohomology.\Box

\proposition {lc and sheaf cohom}
With $P$ and $B$ as in the Theorem,

(a) For $i \geq 1$, there is an isomorphism of graded
$S$-modules
$$\H^i_*(\tilde{P}) \, \cong \, \H^{i+1}_B(P),$$

(b) There is an exact sequence of graded $S$-modules
$$\quad 0 \rTo \H^0_B(P)\rTo P \rTo \H^0_*(\tilde{P}) \rTo
\H^1_B(P) \rTo 0.$$

\proof Both assertions follow from the
fact that
$$
\H^i_B(P)\,\cong\, H^i(C({\bf m}, P)).\Box
$$

This concludes our review of basic machinery.

Because of the similarity of the computation of sheaf cohomology
and local cohomology in this setting, we can translate
the results of the last section directly. For example,
translating \ref{dependence for monomial} we have:

\theorem {dependence on neg p}
The module $\H^i_*(\O_X)$ is $\Z^n$-graded.
If $p, q \in \Z^n$ satisfy $\neg(p) = \neg(q)$, then there is
a canonical isomorphism
$$\H^i_*(\O_X)_p \iso\H^i_*(\O_X)_q.\Box$$

Thus we can define
for each $i \geq 0$,
$$
\Sigma^{[1]}_i=\{ I \subset \{1,\dots,n\} \mid  \H^i_*(\O_X)_{p_I} \neq 0 \},
$$
and we have

\corollary{constancy2}
$$
Supp \H^i_*(\O_X) = \bigcup_{I \in\Sigma^{[1]}_i} L_I.
\Box$$

For a complete toric variety
$\H^i(\O_X(D))$ is finite dimensional
for every invariant divisor $D$ of $X$, so the conditions
of \ref{inter equivs} hold. In the general toric setting
we can provide a more explicit
finiteness condition
to complement those given in
\ref{inter equivs}. We will call a 1-dimensional cone of the
fan $\Delta$ an {\it edge.}

\proposition{inter equivs toric} With notation as in \ref{inter equivs},
if $B$ and $\D$ come from a toric variety $X$ with
fan $\Delta$ as above,
then $\H^i(\O_X(D))$ is finite dimensional
for every invariant divisor $D$ of $X$
iff for every $I\in\Sigma^{[1]}_i$,
and every $p\in\Z^n$, $(p+M)\cap L_I$ is finite.
The conditions (a)-(d) of \ref{inter equivs} are also equivalent
in this case to
\item{(e)} There is no hyperplane $H\subset N_\RR=\RR^d$
such that every edge of $\Delta$ indexed by
an element of $I$
lies in or on one side of $H$ and every
edge
indexed by an element of $I'$
lies in or on the other side of $H$.

\proof The first statement is simply the application
of \ref{inter equivs} to our case.
We continue the notation of \ref{inter equivs} and prove
$(b)\equiv (e)$: If $0\neq p\in M\cap C_I$ is regarded as a
linear functional on $N_R$, then its hyperplane of zeros satisfies
the given condition. Conversely, a hyperplane $H$
that satisfies the condition separates the two cones spanned by
the edges indexed by $I$ and $I'$.
As these cones are both spanned by integral vectors, we may take
$H$ to be integral, and an integral functional $p$ vanishing
on $H$ is in $M\cap C_I$.\Box

In the toric case,
the actual value of
$\H^i_*(\O_X)_p$ can be computed from a topological formula.
Note that Fulton [1993] gives a slightly different topological
description in the special case where the divisor in $X$ corresponding
to $p$ is Cartier.

 If $\Delta$ is the fan of $X$, then
  we write $|\Delta|$
for the union of the cones in $\Delta$.
For a subset $I\subset\{1,\ldots, n\}$, we define $Y_I$ to be the union
of those cones in $\Delta$ having all the edges in the complement of $I$
( if $I=\{1,\ldots, n\}$, we take $Y_I=\{0\}$ ). We use the notation
$\H^i_{Y_I}(|\Delta|) = \H^i(|\Delta|, |\Delta|\setminus Y_I).$ Since
$|\Delta|$ is contractible, we have
$\H^i_{Y_I}(|\Delta|) = \H^{i-1}(|\Delta|\setminus Y_I)$
(reduced cohomology with coefficients in $k$).
We have

\theorem {topological description} With the above notation, if
$p\in\ZZ^n$ and $I={\rm neg}(p)$, then
$$\H^i_*(\O_X)_p\iso\H^i_{Y_I}(|\Delta|).$$

\proof If $i=0$, $\H^0_*(\O_X)=S$, so that $\H^0_*(\O_X)_p\iso k$
if $I=\emptyset$ and $0$ otherwise. Since the same thing holds trivially for
$\H^0_{Y_I}(|\Delta|)$, we can conclude this case.

Suppose now that $i\geq 1$. In this case,
$$\H^i_*(\O_X)_p\iso\H^{i+1}_B(S)_p.$$

We will use the topological description
of $\H^{i+1}_B(S)_p$ for any square-free monomial ideal $B$,
from \mustata, [1999a], Theorem 2.1.
It says that if $m_1,\ldots, m_r$ are the minimal monomial
generators of $B$, $I={\rm neg}(p)$ and $T_I$ is the simplicial complex
on $\{1,\ldots, r\}$ given by
$$T_I=\{J\subset\{1,\ldots, r\} | X_i\not | {\rm lcm}(m_j;j\in J)
{\rm for}\,{\rm some}\, i\in I\},$$
then $$\H^{i+1}_B(S)_p\iso\H^{i-1}(T_I, k).\,\,\,(*)$$
When $I=\emptyset$, $T_I$ is defined as the void complex. Note that in $(*)$
and everywhere below, the cohomology groups are reduced and with
coefficients in $k$.

In our situation the minimal generators of $B$ correspond exactly to the
maximal cones in $\Delta$. Therefore, if $\sigma_1,\ldots,\sigma_r$
are the maximal cones in $\Delta$, then $T_I$ is the simplicial complex
on $\{1,\ldots, r\}$ given by: $J\in T_I$ iff there is $i\in I$
such that $X_i\not |  X_{\hat \sigma_j}$, for every $j\in J$.
Equivalently, $J$ is in $T_I$ iff the intersection of the cones in J
has some edge in $I$.

If $I=\emptyset$, then the conclusion follows trivially
from $(*)$, since the void complex has trivial reduced cohomology.

Therefore we may suppose that $I\neq\emptyset$. We have already seen
that $\H^i_{Y_I}(|\Delta|)\iso\H^{i-1}(|\Delta|\setminus Y_I)$.

Let's consider the following cover of $|\Delta|\setminus Y_I$:
$$\{\sigma\setminus Y_I | \sigma\in\Delta\, {\rm maximal}\,{\rm cone}\}.$$
$(\sigma_{i_1}\setminus Y_I)\cap\ldots\cap(\sigma_{i_k}\setminus Y_I)
\neq\emptyset$ iff $\sigma_{i_1}\cap\ldots\cap\sigma_{i_k}\not\subseteq Y_I$.
This means precisely that $\sigma_{i_1}\cap\ldots\cap\sigma_{i_k}$
has an edge in $I$. Therefore $T_I$ is the nerve of the
above cover.

On the other hand, if $(\sigma_{i_1}\setminus Y_I)\cap\ldots
\cap(\sigma_{i_k}\setminus Y_I)\neq\emptyset$, then this set is
equal to $\sigma\setminus (\sigma\cap Y_I)$, where $\sigma=
\sigma_{i_1}\cap\ldots\cap\sigma_{i_k}\in\Delta$. But a sharp cone
minus some of its proper faces is contractible. Therefore by
Godement [1958], Theorem 5.2.4, we get
$$\H^{i-1}(|\Delta|\setminus Y_I)\iso
 \H^{i-1}(T_I)$$
and finally, using $(*)$,
$$\H_*^i(\O_X)_p\iso \H^i_{Y_I}(|\Delta|).\Box  $$

\example{cohom of surfaces} Let $X$ be a (not necessarily complete)
toric surface. In this case \ref{topological description}
gives a complete description of $\Sigma^{[1]}_i$.

For $i=0$, $\H^0_{Y_I}(|\Delta|) =0$, unless $Y_I=|\Delta|$
i.e. $I=\emptyset$, in which case $\H^0_{Y_I}(|\Delta|)\iso k$.
Therefore $\Sigma^{[1]}_0=\{\emptyset\}$.

For $i=2$, $\H^2_{Y_I}(|\Delta|)\iso\H^1(|\Delta|\setminus Y_I)$,
which is zero, unless $X$ is complete and $I=\{1,\dots, n\}$,
in which case $\H^2_{Y_I}(|\Delta|)\iso k$. Hence $\Sigma^{[1]}_2=
\{\{1,\dots,n\}\}$ if $X$ is complete and it is empty otherwise.

The interesting case is $i=1$, when $\H^1_{Y_I}(|\Delta|)\iso
\H^0(|\Delta|\setminus Y_I)$, so that
$\dim_k\H^1_{Y_I}(|\Delta|)$ is the number of connected components
of $|\Delta|\setminus Y_I$ minus one, unless $Y_I=|\Delta|$, in which
case it is zero.

If $X$ is complete, then $I\not\in\Sigma^{[1]}_1$ iff the
one-dimensional cones in $I$ form a sequence $D_1,\ldots,D_n$
such that $D_i$ and $D_{i+1}$ are adjacent for $1\leq i\leq r-1$.

In particular, we see that it is possible to have some $I$ satisfying
the equivalent conditions in \ref{inter equivs toric}, but such
that $I\not\in\Sigma^{[1]}_i$, for all $i$. Take, for example the
complete toric surface associated to the following fan:

$$\epsfbox{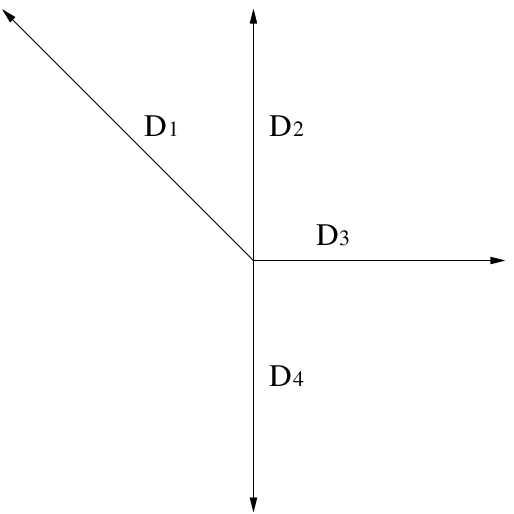}$$

Then $\{1,3,4\}$ has the above property: it satisfies condition (e)
of \ref{inter equivs toric}

\section{cohomology bounds} Bounding $\H^i_B(S)_\delta$ and $\H^i\O_X(D)$

In this section we suppose that $B$ and $D$
are such that the local cohomology
$\H^i_B(S)_\delta$ is finite dimensional for all $\delta\in\D$.
The motivating examples are those when
$B$ and $\D$ come from
the homogeneous coordinate ring of a complete toric variety $X$.
We will derive an effective means
of calculating any $\H^i_B(S)_\delta$ or
$\H^i\O_X(D)$ in terms of an Ext.
The bounds we need come from \ref{empty intersection}.

For simplicity we give the results in this section only for the
local cohomology case. For the corresponding results on
sheaf cohomology one merely replaces $\H^i_B(S)_\delta$
by $\H^i(\O_X(D)$ (where $D$ is a torus-invariant divisor
representing $\delta$) and $\Sigma_i$ by $\Sigma^1_i$.
At the end of the section we illustrate this translation by
working out the example of the two-dimensional rational normal
scrolls (Hirzebruch surfaces).

We adopt
$p_I, C_I$, $\Sigma_i$ and $\Sigma^1_i$ of the previous sections.

\proposition {bound0}
If $I \in \Sigma_i$, and $ p\in\ZZ^n$ then
$$( p + M) \cap (p_I + C_I)$$ is (contained in) a bounded polyhedron.
If we define $f_{I,j}(\delta)$ to be the maximum of zero and the
negative of the minimum
value of the $j$\th\ coordinates
of the points in $( p + M) \cap (p_I + C_I)$, where $p$ is any representative
of $\delta \in\D$, then
$$\Ext^i(S/B^{[\ell]},S)_\delta = H^i_B(S)_\delta$$
if and only if
$\ell \geq \max\{f_{I,j}(\delta) \mid I \in \Sigma_i, j\in I\}$.

\proof
Let $Q = ( p + M) \cap (p_I + C_I)$.
The first statement follows since
$|Q| \leq \dim H^i_B(S)_\delta < \infty$,
so $Q$ consists of the integral points of a convex bounded polyhedron.  The
last
statement follows by \ref{dependence on neg p} and
\ref{lc as union}. \Box

One may find the bound $f_{I,j}(\delta)$ in various ways. For example, one
may solve a linear programming problem in $M \otimes \RR$.  An equivalent
definition of
$f_{I,j}$ is the following: $f_{I,j}(\delta)$ is the least non-negative
integer $\ell$
such that
  $$\delta \not\in \phi(-\ell e_j + L_I).$$
If $\D$ has no torsion, then one may find $f_{I,j}(\delta)$ exactly
by using the facet equations of the cone $\phi(C_I)$.  If $\D$ has torsion,
this method still gives a bound, although possibly not the exact value.
\medskip

\proposition{bound1}
If $C_I \cap M=0$ then the absolute value of each coordinate of
a point of $( p+M)\cap (p_I+C_I)$ is bounded above by
$${d^2 \max_i\{| p_i|\}Q_1Q_{d-1}\over q_d},$$
where
$q_d$ is the minimum of the nonzero absolute values of the $d\times d$
minors of $\rho$ ($q_d\geq 1$, as $\rho$ is a matrix of integers) and
$Q_j$ is the maximum of the absolute values of the
$j\times j$ minors of $\rho$.

\proof Since $p_I+C_I\subset C_I$ it suffices to bound the points
of $( p+M)\cap C_I$. Since this is a bounded polyhedron,
the coordinates of its points, as vectors of $M$, are bounded by the
coordinates
of its vertices. Each vertex $u\in M$ is defined by a system of $d$
linear equations in the $d$ coordinates of $u$
having the form $\rho'u  = - p$,
where $\rho'$ is a $d\times d$ submatrix of $\rho$ with nonvanishing
determinant. Thus the $i^\th$ coordinate $u_i$ has the form
$$
{\sum_{j=1}^d -{ p}_j a_j\over{\rm det}\rho'}
$$
where the $a_j$ are $(d-1)\times(d-1)$ minors of $\rho'$, and
thus has absolute value bounded by
$(d \max_i\{| p_i|\}Q_{d-1})/q_d.$
Applying $\rho$, we get a vector in $\ZZ^n$ whose coordinates
have absolute value at most $dQ_1$ times this, whence the
given bound.
\Box

\corollary{computation1} Suppose $p\in\ZZ^n$ is
a representative for $\delta\in \Ch$.
If
$$\ell = {d^2 \max_i\{| p_i|\}Q_1Q_{d-1}\over q_d}$$ is
the bound in \ref{bound1} then
$$
\H^i_B(S)_\delta = \Ext^i_S(S/B^{[\ell]}, S)_{\delta}.
$$

\proof Combine \ref{bound0} and \ref{bound1}. \Box

\corollary{computation2} Fix an integer $i \geq 0$, and an element $\delta
\in \Ch$.
The inclusion
  $$Ext^i(S/B^{[\ell]},S)_\delta \longrightarrow H^i_B(S)_\delta$$
is an isomorphism if and only if
  $$\ell \geq f_{I,j}(\delta),$$
for all $I \subset \Sigma_i$ and all $j \in I$
where $f_{I,j}(\delta)$ is the bound defined in \ref{bound0}.
\bigskip

The same result may be looked at from a different point of view:

\corollary{} Fix an integer $i \geq 0$, and an integer $\ell \geq 0$.
The set
of $\delta \in \Ch$ such that
$$
Ext^i(S/B^{[\ell]},S)_\delta \longrightarrow \H^i_B(S)_\delta
$$
is an
isomorphism is the complement (in $\Ch$) of the union of the
translated pointed cones
  $$\phi( - \ell e_j) + \phi(L_I),$$
for all $I \in \Sigma_i$ and all $j \in I$.
\Box

\example{scroll} {\bf [Surface scrolls]}
Choose a non-negative integer $e \geq 0$.  Let
  $$0 \rTo \Z^2 \rTo^\rho_{\pmatrix{
          1 & 0 \cr
          0 & 1 \cr
          -1 & e \cr
          0 & -1\cr}}
    \Z^4 \rTo^\phi_{\pmatrix{
          0 & 1 & 0 & 1\cr
          1 & 0 & 1 & e\cr}}
    \Ch = \Z^2 \rTo 0.$$
The matrix $\rho$ defines a fan $\Delta \subset {\bf R}^2$, with $X =
X(\Delta)$
a rational normal surface scroll.  The ring $S = k[x_1, x_2, x_3, x_4]$ has the
$\Z^2$-grading $\deg(x_1) = \deg(x_3) = (0,1)$, $\deg x_2 = (1,0)$, and
$\deg x_4 = (1,e)$.  The irrelevant ideal $B = (x_1, x_3) \cap (x_2, x_4)$.
The non-zero Ext modules of $S/B$ (shown with their fine gradings) are
  $$\Ext^2(S/B,S) = S(-1,0,-1,0)/(x_1, x_3) \oplus S(0,-1,0,-1)/(x_2,x_4),$$
and
  $$\Ext^3(S/B,S) = S(-1,-1,-1,-1)/(x_1,x_2,x_3,x_4).$$
Therefore the non-empty sets $\Sigma_i$ are
  $$\Sigma^{[1]}_1  = \Sigma_2 = \{ \{1,3\}, \{2,4\} \},$$
and
  $$\Sigma^{[1]}_2  = \Sigma_3 = \{ \{1,2,3,4\} \}.$$
To compute the cohomology module $\H^i(\O_X(a,b))$, one may use
$ p = (b,a,0,0) \in \Z^4,$ and the bound in \ref{computation1} becomes
  $$\ell = 4e^2 \max\{|a|,|b|\}.$$
This is in general not best possible.  For $i = 0$, one may always take $\ell =
0$.  For $i = 2$, by solving for the minimum of any coordinate of any element
of $( p + M) \cap (p_I + C_I)$, where $I = \{1,2,3,4\}$, one finds that
  $$\ell = \max(-a-1, -b-e-1,0)$$
is the best possible bound.

For $i=1$, we have two sets $I$ to consider.  For $I = \{1,3\}$, the set of
$(a,b)$ for which $( p + M) \cap (P_I + C_I)$ is non-empty is the set
$$\{ (a,b) \mid a \geq 0, b \leq ae-2\}.$$
The corresponding bound is the minimum of coordinates 1,3 of any point in
$(p + M) \cap (P_I + C_I)$.  This works out to a bound of
  $$\ell = \max(-b+ae-1, 0),$$
for any $(a,b)$ in the above set.
The second set, $I = \{2,4\}$, gives rise to a region
$$\{ (a,b) \mid a \leq -2, b \geq ae+e\}.$$
The corresponding bound for points in this region is
  $$\ell = \max(-a-1, b/e-a, 0).$$

For the case $e=1$,
The illustration below shows
 the cones $\phi(L_I)$, for $I =
{1,3}, {2,4},
$ and ${1,2,3,4}$.  In addition, for each of these three translated cones, the
lines $\ell = r$ are shown, for $r = 2,3,4,5$.  For $\ell = 1$, the only points
in these three cones such that $\Ext^i(S/B^{[\ell]},S)_\delta =
\H^i_B(S)_\delta$
are the vertices of the cones.

\epsfbox{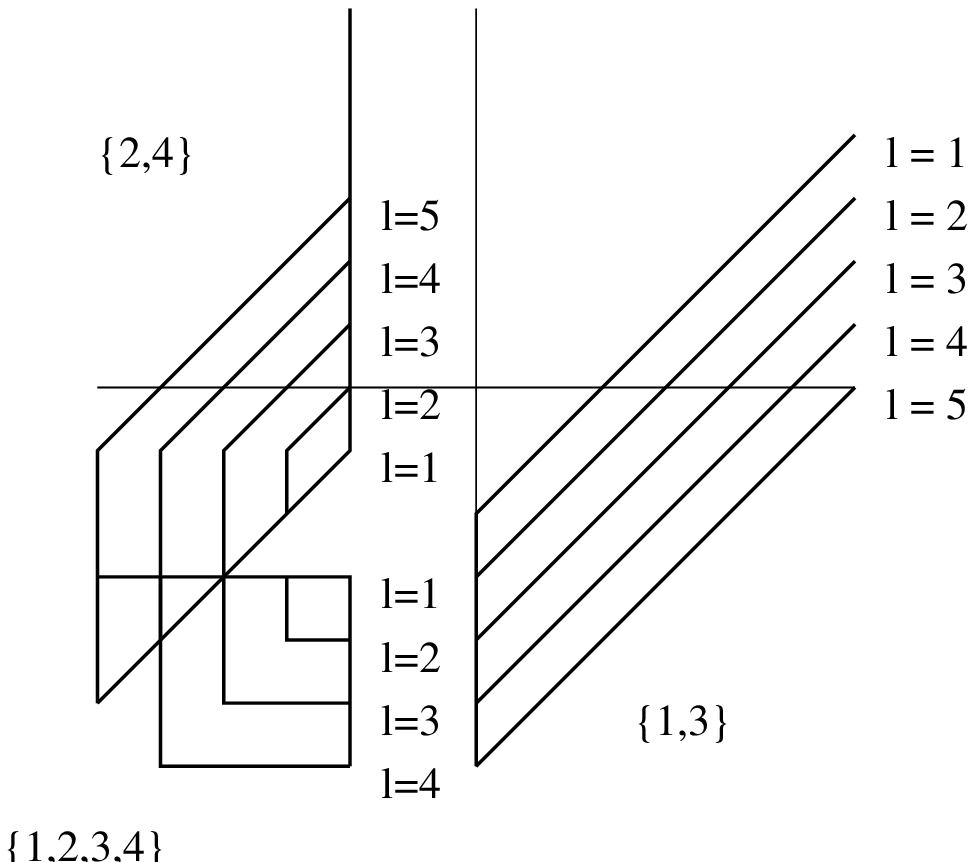}

\section{modulebounds}  Bounds for $\H^i_B(P)_\delta$ and
$\H^i(\tilde{P}(\delta))$

Once we know bounds for computing the local cohomology
of the ring $S$, we can deduce the bounds for the case of an arbitrary
finitely generated graded module. This is done in the proposition below
using a spectral sequence argument which has the same flavor
as the arguments in Smith [1999].

\smallskip

\proposition{bounds-by-functions} Suppose that we have functions $f_i:
\Ch\longrightarrow {\bf N}^*$ such that:
$$\Ext^i_S(S/B^{[k]}, S)_\delta\longrightarrow \H_B^i(S)_\delta$$
is an isomorphism for every $\delta\in\Ch(X)$, $i\geq 0$ and $k\geq
f_i(\delta)$.
Let $P$ be a finitely generated $\D$-graded $S$ module.
Let $F_{\bullet}$ be a free resolution of $P$ with
$$F_i=\bigoplus_{\alpha\in\Ch(X)}S(-\alpha)^{\beta_{i,\alpha}}.$$
Then
$$\Ext^i_S(S/B^{[k]}, P)_\delta\longrightarrow \H^i_B(P)_\delta$$
is an isomorphism if for every $j\geq 0$ and every $\alpha$ such that
$\beta_{j,\alpha}\neq 0$ we have
$$k\geq\max\{f_{i+j-1}(\delta-\alpha), f_{i+j}(\delta-\alpha)\}$$
if $j\geq 1$ and $k\geq f_i(\delta-\alpha)$ if $j=0$.

\proof Let $E^{\bullet}$ be the complex given by $E^{-i}=F_i$, $i\geq 0$.
We will use the two spectral sequences of hypercohomology for the
functors $\Hom(S/B^{[k]},-)$ and $\H_B^0(-)$ and the complex
$E^{\bullet}$ (see Weibel [1994] for details about hypercohomology).

For $\Hom(S/B^{[k]},-)$ we get the spectral sequences:
$$'E^{p,q}_2=\H^p(\Ext^q_S(S/B^{[k]}, E^{\bullet}))\Rightarrow
{\bf Ext}^{p+q}(S/B^{[k]}, E^{\bullet})$$
$$''E^{p,q}_2=\Ext^p_S(S/B^{[k]}, \H^q(E^{\bullet}))\Rightarrow
{\bf Ext}^{p+q}_S(S/B^{[k]}, E^{\bullet}).$$

Since $\H^q(E^{\bullet})=0$ for $q\neq 0$ and $\H^0(E^{\bullet})=P$,
the second spectral sequence collapses and the first one becomes:
$$'E^{p,q}_2=\H^p(\Ext^q_S(S/B^{[k]}, E^{\bullet}))\Rightarrow
\Ext^{p+q}_S(S/B^{[k]}, P).$$

Applying the same argument for $H^0_B(-)$, we get the spectral sequence:
$$'\overline{E}^{p,q}_2=\H^p(\H^q_B(E^{\bullet}))\Rightarrow
\H^{p+q}_B(P).$$

Moreover , the natural map $\Hom(S/B^{[k]},-)\longrightarrow
H^0_B(-)$ induces a natural map of spectral sequences.
Therefore, in order for the map
$$\Ext^i_S(S/B^{[k]}, P)_\delta\longrightarrow\H^i_B(P)_\delta$$
to be an isomorphism, it is enough to have
$$u^{p,q} : \H^p(\Ext^q_S(S/B^{[k]}, E^{\bullet}))_\delta\longrightarrow
\H^p(H^q_B(E^{\bullet}))_\delta$$
isomorphism for every $p$ and $q$ with $p\leq 0$ and $p+q=i$.

In the diagram below:
$$\diagram[midshaft]
\Ext^q_S(S/B^{[k]}, E^{p-1})_\delta&\rTo&\Ext^q_S(S/B^{[k]}, E^p)_\delta&\rTo&
\Ext^q_S(S/B^{[k]}, E^{p+1})_\delta& \\
\dTo v^{p-1,q}& &\dTo v^{p,q}& &\dTo v^{p+1,q}& \\
\\H^q_B(E^{p-1})_\delta&\rTo&\H^q_B(E^p)_\delta&\rTo&\H^q_B(E^{p+q})_\delta&\\
\enddiagram $$
the vertical maps are injective by Theorem 0.2. Therefore,
$u^{p,q}$ is an isomorphism if both $v^{p,q}$ and $v^{p-1,q}$ are
isomorphisms.

Since $E^p=\oplus_{\alpha}S(-\alpha)^{\beta_{-p,\alpha}}$, using
the properties of the functions $f_i$ we get the conclusion of the proposition.
\Box

Translating this to sheaf cohomology is immediate:

\corollary{bounds-by-functions2} Suppose that we have functions
$f_i: \Ch\longrightarrow {\bf N}^*$
such that:
$$\Ext^i_S(B^{[\ell]}, S)_\delta\longrightarrow \H^i(\O_X(\delta))$$
is an isomorphism for every $\delta\in\Ch$, $i\geq 0$ and
$\ell\geq f_i(\delta)$.
Let $P$ be a finitely generated graded $S$ module.
Let $F_{\bullet}$ be a minimal free resolution of $P$ with
$$F_i=\bigoplus_{\alpha\in\Ch}S(-\alpha)^{\beta_{i,\alpha}}.$$
Then
$$\Ext^i_S(B^{[\ell]}, P)_\delta\longrightarrow H^i_*(\tilde{P})_\delta$$
is an isomorphism if for every $j\geq 0$ and every $\alpha$ such that
$\beta_{j,\alpha}\neq 0$ we have
$$k\geq\max\{f_{i+j-1}(\delta-\alpha), f_{i+j}(\delta-\alpha)\}$$
if $j\geq 1$ and $k\geq f_i(\delta-\alpha)$ if $j=0$.\Box

\section{alg} {\bf Algorithms and Examples}

In this section we give explicit algorithms corresponding to the
theorems in the previous section, and then give some examples.
These algorithms have been implemented in the {\it Macaulay 2\/}
system (Grayson-Stillman [1993--]).

Given the ideal $B$ and the grading $\phi : \Z^n \rTo \D$,
we once and for all compute $\Sigma_i = \Sigma_i(B)$, the $\Z^n$-support
of $\Ext^i_S(B,S)$, as well as the
bounding hyperplanes for the cones $\phi(C_I)$, for each $I \in
\Sigma_i$.

The following algorithm is essentially Corollary 3.1 in \mustata [1999a].

\algorithm {gbalg1} {\bf [$\Sigma_i(B)$, all $i$]}
\cleartabs
\+ {\bf input: } A square-free
monomial ideal $B \subset S$.\cr
\+ {\bf output:} The $\ZZ^n$-support of the module $\Ext^i(B,S)$,
for each $i$.
\cr
\+ {\bf begin}\cr
\+ \qquad & {\bf set} $B^*$ := Alexander dual of $B$:\cr
\+ & \qquad & $B^*$ := the monomial ideal generated by the
square-free monomials in $(x_1^2, \ldots, x_n^2) : B$.  \cr
\+ & {\bf set} $\Sigma_i$ := $\emptyset$, for all $i$\cr
\+ & {\bf for} $i := 0$ {\bf to} $n$ {\bf do}\cr
\+ & \qquad & {\bf set} $T := Tor_i(S/B^*,k)$ \cr
\+ & & {\bf for each} $\ZZ^n$ degree $p$ of $T$ {\bf do}\cr
\+ & & \qquad & {\bf set} $I := \{ j \in 1..n \mid p_j \neq 0 \}$\cr
\+ & & & $\Sigma_{(|p|-i+1)}$ := $\Sigma_{(|p|-i+1)} \cup
\{I\}$\cr
\+ & {\bf return} the sets $\Sigma_i$, for $i=1,\ldots,n$\cr
\+ {\bf end.}\cr
\medskip

The following routine corresponds to \ref{computation2}.  The ideal $B$ and
grading $\phi$ are implicit parameters.

\algorithm {gbalg2} {\bf [bound$(i,\delta,S)$]}
\cleartabs
\+ {\bf input: } An integer $i \in \ZZ$, and a coarse degree $\delta
\in \Ch$.\cr
\+ {\bf output:} The least integer $\ell$ such that
$\Ext^i_S(S/B^{[\ell]}, S)$ in degree $\delta$ is equal to
$\H^i_B(S)_\delta$.
\cr
\+ {\bf begin}\cr
\+ \qquad & {\bf set} $\ell$ := 0\cr
\+ & {\bf for each} $I$ {\bf in} $\Sigma_i(B)$ {\bf do}\cr
\+ & \qquad & {\bf for each} $j$ {\bf in} $I$ {\bf do}\cr
\+ & & \qquad & Solve the linear programming problem:\cr
\+ & & & \qquad & $m$ :=  the negative of the minimum value of the $j$\th\
coordinate
in $(p + M) \cap L_I$,\cr
\+ & & & & where $p \in \Z^n$ is a representative of $\delta$.\cr
\+ & & & {\bf if} $m > \ell$ {\bf then} {\bf set} $\ell$ := $m$.\cr
\+ & {\bf return} $\ell$ \cr
\+ {\bf end.}\cr
\medskip

In view of the remarks after \ref{bound0}, we may
compute the bound
$m$ in the above algorithm by instead computing the facet equations
of $\phi(C_I)$,
and then using them to find the minimum $m \geq 0$ such that
$\delta \not\in \phi(p_I - m e_j + C_I)$.

\algorithm {gbalg3} {\bf [bound$(i,\delta,F)$]}
\cleartabs
\+ {\bf input: } An integer $i \in \ZZ$, and a coarse degree $\delta
\in \Ch$, and a graded free $S$-module\cr
\+ \qquad $\displaystyle{
{F=\bigoplus_{\alpha\in\Ch}S(-\alpha)^{\beta_{\alpha}}}}$\cr
\+ {\bf output:} The least integer $\ell$ such that
$\Ext^i_S(S/B^{[\ell]}, F)$ in degree $\delta$ is equal to
$\H^i_B(F)_\delta$.
\cr
\+ {\bf begin}\cr
\+ \qquad & {\bf return} the maximum of the numbers\cr
\+ & \qquad & $\{ \hbox{bound}(i,\delta-\alpha,S) \mid \beta_\alpha \neq
0\}$\cr
\+ {\bf end.}\cr
\medskip

\algorithm {gbalg4} {\bf [bound$(i,\delta,P)$]}
\cleartabs
\+ {\bf input: } An integer $i \in \ZZ$, a coarse degree $\delta
\in \Ch$, and a graded  $S$-module $P$.\cr
\+ {\bf output:} An integer $\ell$ such that
$\Ext^i(S/B^{[\ell]}, P)$ in degree $\delta$ is equal to
$\H^i_B(P)_\delta$.
\cr
\+ {\bf begin}\cr
\+ \qquad & Compute a minimal graded free resolution
$F_{\bullet}$ of $P$.\cr
\+ & {\bf set} $\ell$ := bound$(i,\delta,F_0)$\cr
\+ & {\bf for each} $j \geq 1$ {\bf do}\cr
\+ & \qquad & {\bf set} $\ell$ :=  max($\ell$,
bound$(i+j,\delta,F_j)$\cr
\+ && {\bf set} $\ell$ :=  max($\ell$,
bound$(i+j-1,\delta,F_j)$\cr
\+ & {\bf return} $\ell$\cr
\+ {\bf end.}\cr
\medskip

For a specific module $P$, this bound is not always best
possible.  Once we have this bound $\ell$,
we compute $\Ext^i(S/B^{[\ell]}, P)$, using standard methods.  The degree
$\delta$
part
of this module is  easily extracted using a Gr\"obner basis of this module, or
the dimension may be found by computing its (multi-graded) Hilbert
function.

\example{multi-proj1} Let $Y$ be $\P(3,3,3,1,1,1)$ and $X$ its
desingularisation.
Then the homogeneous coordinate ring of $X$ is
$S=k[X_1,\ldots,X_6,T]$, $\Ch(X)=\Z\times\Z$ and ${\rm deg}(X_i)=
(3,1)$ for $1\leq i\leq 3$, ${\rm deg}(X_i)=(1,0)$ for $4\leq i\leq 6$
and ${\rm deg}(T)=(0,1)$.

Let $\Delta_Y$ be the fan defining $Y$. Then the maximal cones
of $\Delta$ are $\sigma_1,\ldots,\sigma_6$, where $\sigma_i$ is generated
by $\overline e_1,\ldots,\hat{\overline e_i},\ldots,\overline e_6$
($\overline e_1+\ldots+\overline e_6=0$ and the lattice $N$ is
generated by $1/3\overline e_1$, $1/3\overline e_2$, $1/3\overline e_3$,
$\overline e_4$, $\overline e_5$ and $\overline e_6$).
Let $\overline f=1/3(\overline e_4+\overline e_5+\overline e_6)$
and for each $i$, $1\leq i\leq 3$ we consider
$\sigma_i=\sigma_{i4}\cup\sigma_{i5}\cup\sigma_{i6}$, where
$\sigma_{ij}$ is obtained by replacing $\overline e_j$ with
$\overline f$ in $\sigma_i$.
Then the maximal cones of the fan $\Delta$ defining $X$ are
$\sigma_{ij}$ for $1\leq i\leq 3$, $4\leq j\leq 6$ and $\sigma_i$
for $4\leq i\leq 6$.

Note that in $\Delta$ any two edges are contained in a maximal cone,
while the only three edges that do not belong to a maximal cone are
$\{\overline e_4,\overline e_5,\overline e_6\}$. Using this and our
topological description of the support we deduce that
$\Sigma_i=\emptyset$ if $i\neq 3$, $4$ or $6$, while
$$\Sigma_3=\{\{\overline e_4,\overline e_5,\overline e_6\}\},$$
$$\Sigma_4=\{\{\overline e_1,\overline e_2,\overline e_3,\overline f\}\},$$
$$\Sigma_6=\{\{\overline e_1,\ldots,\overline e_6,\overline f\}\}.$$

We obtain the functions $f_i:\Ch(X)\longrightarrow{\bf N}^*$ that satisfy
the property in the above Corollary 4.2:
$f_i\equiv 1$ if $i\neq 3$, $4$ or $6$, and
$$f_3(\delta_1,\delta_2)=\cases{1, & if $\delta_2<{\rm max}
                                 (0, 1/3\delta_1+1)$;\cr
                                3\delta_2-\delta_1-2, & if $\delta_2
                                \geq {\rm max}(0,
                                  1/3\delta_1+1)$\cr},$$
$$f_4(\delta_1,\delta_2)=\cases{1, & if $\delta_2+1>{\rm min}
                                  (-3, 1/3\delta_1)$;\cr
                                -\delta_2-3, & if $\delta_2+1\leq {\rm min}
                                  (-3, 1/3\delta_1)$\cr},$$
$$f_6(\delta_1,\delta_2)=\cases{1, & if $\delta_1>-12$ or $\delta_2>-4$;\cr
                                {\rm max}(-\delta_2-3, -\delta_1-11), &
                             if $\delta_1\leq-12$ and $\delta_2\leq-4$\cr}.$$

%%%%%%%%%%%%%%%%%%%%%%%%%%%%%%%%%%%%%%%%%%%%%%%%%%%%%%%
\section {problem on lc} A Problem on Local Cohomology

Let $S$ be a polynomial ring and let $I\subset S$
be an ideal. The natural maps $\Ext^j_S(S/I^d,S)\to\H^j_I(S)$ are
rarely injections (they are almost never surjections). It is thus
interesting and, from the point of view of computation,
useful to ask when some analogue of \mustata's \ref{lc as union}
holds:

\question{Q1} For which ideals $I$
does there exist a sequence of ideals
$$
I\supset I_1\supset I_2\supset\dots\supset I_d\supset\dots
$$
such that each $\Ext^j_S(S/I_d,S)$ injects into $\H^j_I(S)$
and
$$\H^j_I(S)=\bigcup \Ext^j_S(S/I_d,S)?$$

An even more basic question is:

\question{Q2} For which ideals $I\subset S$ is the natural
map $\Ext^j(S/I,S)\to\H^j_I(S)$ an inclusion?

By \ref{lc as union}, Frobenius powers of reduced monomial
ideals do satisfy the condition of \ref{Q2}, and thus
sequences of ideals as in \ref{Q1} exist whenever $I$ is
monomial, but
we do not know the
full answer to \ref{Q2} even for monomial ideals.

\example{unmixed} If we require the condition in \ref{Q2} for all $j$, then
$I$ must be unmixed.
(Proof:
Take $j$ equal to the codimension of the embedded prime.
Both Ext and local cohomology
localize, so we may begin by localizing
and assume that the embedded prime is the maximal ideal.
We then get
$\Ext^j(S/I,S)\neq 0$ but
$\H^j_I(S)=0$ as
one sees by writing
$\H^j_I(S) = \lim
 \Ext^j_S(S/I^{(\ell)},S)$,
where $I^{(\ell)}$ denotes the intersection of the $\ell^\th$
symbolic powers of the minimal primes of $I$.)
But even when
$S$ is a polynomial ring and $I$ is a monomial ideal, this
condition does not suffice. If $I\subset k[a,b,c,d]$ denotes
the monomial ideal $(ab,acd,bd^2,cd^2)=(b,d)\cap(b,c)\cap(a,d^2)$
then one sees by computation that the map
$\Ext^3_S(S/I,S)\to \Ext^3_S(S/I^2,S)$ has kernel $k$.

Craig Huneke pointed out to us one way that one might
produce ideals $I$ for which no sequence as in \ref{Q1}
can exist:
Let $k$ be a
field of characteristic $p>0$, and
let $I$ be a prime ideal in $S$
such that
$\H_I^j(S) = 0$ for $j>\codim(I)$.
Let $J$ be any ideal with
radical $I$. It follows that $\Ext^*(S/J,S)$ injects in
$\H^*_I(S)=\H^*_J(S)$
iff $J$ is perfect---that is
$S/J$ is Cohen-Macaulay (in this case
we say that $I$ is set-theoretically Cohen-Macaulay).
If such a $J$ exists,
then the ideal $J^{[p^n]}$
generated by the ${p^n}^\th$ powers of elements of $J$ is perfect
too (apply the Frobenius to the whole free resolution of $J$,
using the characterization that $J$ is perfect iff it has a
free resolution of length equal to $\codim(J)$),
so the ideals $J_d=J^{[p]}$ satisfy the condition of \ref{Q1}
above.

Suppose now that $S/I$ is $F$-pure. We
claim that if there exists a perfect ideal $I$ with the
same radical as $J$, then $I$ would be perfect. To see this,
suppose $x,y$ is a system of parameters modulo $I$,
and $r,s\in S$ are such that $rx+sy\in I$. It follows that
$$
r^{p^m}x^{p^m}+s^{p^m}y^{p^m}\in J
$$
for large $m$, and thus $s^{p^m}\in (x^{p^m}+J)\subset (x^{p^m}+I)$.
 From $F$-purity we get $s\in (x+I)$, so $x,y$ is a regular sequence.

Thus if $I$ is an imperfect prime ideal such that
$S/I$ is $F$-pure and $H^j_I(S)=0$ for all $j>\codim(I)$
then no sequence as in \ref{Q1} can exist for $I$.
Unfortunately we do not know whether such a prime ideal
exists. See Huneke and Lyubeznik [1990] for some results
where at least the cohomology vanishing is proven.
%\vfill\eject
\bigskip\bigskip

\centerline {\bf References}
\baselineskip=12pt
\parindent=0pt
\frenchspacing
\medskip
%Models:
%\item{} V.~I.~Arnold: A-graded algebras and continued fractions,
%{\sl Communications in Pure and Appl. Math.} {\bf 42} (1989) 993-1000.
%\medskip
%\item{} D.~Cox, J.~Little, D.~O'Shea:
%{\sl Ideals, Varieties and Algorithms},
%Springer, New York, 1992.
%\medskip

%\references

\item{} D.~Bayer, M.~Stillman:
Macaulay: A system for computation in
        algebraic geometry and commutative algebra
Source and object code available for Unix and Macintosh
        computers. Contact the authors, or download from
{\tt ftp:// math.harvard.edu/macaulay} via anonymous ftp.
\medskip

\item{} D.~Cox:
{\sl The homogeneous coordinate ring of a toric variety},
J.Algebraic Geom. 4 (1995), 17--50.
\medskip
\item{} M.~Brodmann and R.~Sharp.
{\sl Local cohomology: an algebraic introduction with geometric
applications}, Cambridge University Press, 1998.
\medskip

\item{} D.~Eisenbud:
{\sl Commutative Algebra with a View Toward Algebraic Geometry},
Springer, New York, 1995.
\medskip

\item{} W.~Fulton: {\sl Introduction to Toric Varieties},
  Annals of Mathematical Studies 131, Princeton University Press, 1993.
\medskip

\item{} R.~Godement:
{\sl Topologie Algebrique et Theorie des Faisceax},
Hermann, Paris, 1958.
\medskip

\item{} D.~Grayson, M.~Stillman:
{\sl Macaulay 2: A system for computation in
        algebraic geometry and commutative algebra}.
Source and object code available at
{\tt http:// www.math.uiuc.edu/Macaulay2} (1993--).
\medskip

\item {} A.~Grothendieck: {\sl Local Cohomology}, Springer Lecture Notes in
Math
41, Springer-Verlag, Heidelberg, 1967.
\medskip

\item {} C.~Huneke, G.~Lyubeznik:
{\sl On the vanishing of local cohomology modules},
Invent.Math. 102 (1990) no.1, 73--93.
\medskip

\item{} M.~\mustata :
{\sl Local Cohomology at Monomial Ideals}, J. Symbolic Computation,
to appear, 1999a.
\medskip

\item{} M.~\mustata :
{\sl Vanishing Theorems on Toric Varieties}, preprint, 1999b.
\medskip

\item{} G.~Smith: {\sl Computing global extension modules for coherent
sheaves on a projective scheme}, J. Symbolic Computation, to appear, 1999.
\medskip

\item{} W.~Vasconcelos:
{\sl Computational Methods in Commutative Algebra and Algebraic
Geometry}, Algorithms and computation in mathematics 2,
Springer-Verlag, 1998.
\medskip

\item{} C.~Weibel:
{\sl An Introduction to Homological Algebra}, Cambridge
        studies in advanced mathematics 38, Cambridge University Press,
        1994.

\bigskip
\vbox{\noindent Author Addresses:}
\smallskip
\noindent{David Eisenbud}\par
\noindent{Department of Mathematics, University of California,
Berkeley CA 94720}\par
\noindent{eisenbud@math.berkeley.edu}
\smallskip
\noindent{Mircea \mustata}\par
\noindent{Department of Mathematics, University of California,
Berkeley CA 94720}\par
\noindent{mustata@math.berkeley.edu}
\smallskip
\noindent{Mike Stillman}\par
\noindent{Department of Mathematics, Cornell University
Ithaca, NY 14853-0001}\par
\noindent{mike@math.cornell.edu}\par

\end

%%%%%%%%%%%%%%%%jetsam%%%%%%%%%%%%%%%%%%%%%%%%%%%%%%%%%%%%%%%%%%%%

Given any subset $I\subset \{1,\dots,n\}$, let $p_I \in \ZZ^n$ be the vector
with
$$
(p_I)_j = \cases{-1, &if $j \in I$\cr
                        0,  &if $j \not\in I$}
$$

\section {surface case} The Surface Case

\fix{Redo this in terms of the bound above and the cohom
calculation in bound2}

Let $\Delta \subset {\bf R}^2$ be an integral fan in the plane, having $n$
edges.
%As above, let
%$$\ZZ^2 {A \above \longrightarrow} \ZZ^n {\phi \above \longrightarrow} D =
%coker A.$$

Our first goal is to describe the ``bad'' locus for the irrelevant ideal $B$
of this fan.  We do this by explicit computation using $B$.

Let $R = k[x_1, \ldots, x_n]$.  Indices are computed modulo $n$,
e.g. $x_0 = x_n$ and $x_{n+1} = x_1$. Define $B = (g_1, \ldots, g_n)$,
where
$$g_i = {\prod_{j=1}^n x_j \over x_i x_{i-1}}.$$

Let $\ZZ^n$ be generated by the unit vectors $e_1, \ldots, e_n$.  With this
notation, $\deg x_i = e_i$.  Indices of the $e_i$ are computed modulo $n$
as well.

\proclaim Lemma 1. $R/B$ has the (fine-graded) resolution
$$\matrix{0 & \mapright{} & F_3 & \mapright{\partial_3} & F_2
&\mapright{\partial_2} &F_1
&\mapright{\partial_1} &R& \mapright{} &R/B &\mapright{} &0},$$
where
$$\eqalign{
F_1 &= \bigoplus_{i=1}^n R(-e_1 - \cdots - e_n + e_i + e_{i-1}),\cr
F_2 &= \bigoplus_{i=1}^n R(-e_1 - \cdots - e_n + e_i),\cr
F_3 &= R(-e_1 - e_2 - \cdots -e_n),\cr}$$
and the maps are
$$\eqalign{
\partial_1 &= \pmatrix{g_1 & g_2 & \ldots & g_n},\cr
\partial_2 &= \pmatrix{x_2 & -x_1 \cr
                        & x_3 & -x_2 \cr
                        &     & \ddots \cr
                        &     &     & x_n & -x_{n-1} \cr
                    -x_n &    &     &     & x_1 \cr},\cr
\partial_3 &= \pmatrix{x_1 \cr x_2 \cr \vdots \cr x_n \cr}\cr.}$$

\proof\ Direct computation. \Box

The resolution of $R/B^{[\ell]}$ is obtained by replacing each $x_i$ with
$x_i^\ell$ and each $e_i$ with $\ell e_i$.

For $1 \leq i < j \leq n$, we say that the pair $\{i,j\}$ is {\sl adjacent}
if $j = i+1$, or $j = n, i=1$.  Otherwise, we call the pair $\{i,j\}$ {\sl
non-adjacent}.

\proclaim Lemma 2. $\Ext^1(B,R)$ is generated by ${n \choose 2} - n$
elements
$$f_{i,j} = x_j a_i^* - x_i a_j^*,$$
where $\{i,j\}$ is non-adjacent pair, and $a_1^*, \ldots, a_n^*$ is the dual
basis to the given basis of $F_2$.  The degree of $f_{i,j}$ is
$$\deg f_{i,j} = e_i + e_j - \sum_{m=1}^n e_m.$$

\proof\ The kernel of $\partial_3^*$ is generated by the Koszul syzygies on
the variables.  The image of $\partial_2^*$ consists of those Koszul syzygies
corresponding to adjacent pairs of variables.  The remaining Koszul
syzygies generate,
and these are precisely the $f_{i,j}$ above.  A simple computation yields
their fine degrees. \Box

Similarly, the generators of $\Ext^1(B^{[\ell]},R)$, and their degrees, are
obtained by replacing
each $x_i$ with $x_i^\ell$, and each $e_i$ with $\ell e_i$.

For a non-adjacent pair $\{i,j\}$, define the ideal $P_{i,j}$ to be
$$P_{i,j} = (x_{i+1} \ldots x_{j-1},\ x_{j+1} \ldots x_n x_1\ldots x_{i-1}).$$
The components of this ideal are all codimension two, and occur as minimal
primes of $B$.

\proclaim Lemma 3. Let $1 \leq i < j \leq n$ be a non-adjacent pair, and let
$f_{i,j}$ be the generator of $\Ext^1(B,R)$ above.  Then
$$Ann(f_{i,j}) = P_{i,j}.$$

\proof\ The inclusion $B \subset P_{i,j}$ gives rise to the following portion
of the long exact sequence for Ext:
$$\Ext^1(P_{i,j}/B,R) \rightarrow \Ext^1(P_{i,j}, R) \rightarrow \Ext^1(B,
R),$$
Since $P_{i,j}/B$ has codimension two, the first term is zero, and so
there is an inclusion
$$\Ext^1(P_{i,j}, R) \rightarrow \Ext^1(B,R).$$
However, $\Ext^1(P_{i,j}, R) \cong R/P_{i,j}$, and its generator maps to
$f_{i,j}$.
\Box

Similarly, the annihilator of the corresponding generator $f_{i,j}^{[\ell]}$ of
$\Ext^1(B^{[\ell]},R)$ is $P_{i,j}^{[\ell]}$.

We now describe the support of the $R$-module $H^1(\tilde{R})$.  The support
of this module is the union of the supports of
$\Ext^1(B^{[\ell]},R)$, over all $\ell$.  Let $\Sigma$ denote the set of
subsets $I$ of $\{1, 2, \ldots, n\}$ whose complement contains a
non-adjacent pair.
Given a subset $I \in \Sigma$, let $p_I \in \ZZ^n$ be the vector
with $$(p_I)_i = \cases{-1, &if $i \in I$\cr
                        0,  &if $i \not\in I$.\cr}$$
Define $C_I \subset \ZZ^n$ to be the additive monoid generated by
$\epsilon_I(i) e_i$, for $1 \leq i \leq n$, where
$$\epsilon_I(i) = \cases{-1, &if $i \in I$\cr
                        1,  &if $i \not\in I$.\cr}$$

\proclaim Lemma 4. With this notation,
$$\bigcup_{\ell = 1}^\infty Supp(\Ext^1(B^{[\ell]},R)) = \bigcup_{I \in
\Sigma} p_I + C_I.$$

\proof\ This proof needs to be added! \Box

\proclaim Corollary. The bad locus is
$$Bad_{i,\ell}(B) = \bigcup_{I \in \Sigma} (p_I + C_I) \bigcap \{ p \in
\ZZ^n \mid some\ component\ of\ p \ is\ < -\ell\}.$$

\proclaim Corollary. The bad locus is
$$Bad_{i,\ell}(B) = \bigcup_{I \in \Sigma, i \in I} - \ell e_i + p_I + C_I.$$

\proof\ Follows immediately from the previous section and Lemma 4. \Box

\bigskip

%%From section 1
We wish to understand for which fine degrees $p$ the natural map
$$\alpha_p : Ext^i(B^{[\ell]},R)_p \longrightarrow H^i(\tilde{R})_p$$
is an isomorphism.  To this end, define the ``bad'' locus to be
$$Bad_{i,\ell}(B) := \{ p \in {\bf Z}^n \mid \alpha_p {\rm is\ not\ an\
isomorphism}\}.$$

\proclaim Corollary 2. The ``bad'' locus of $B$ is
$$Bad_{i,\ell}(B) = \{ p=(p_1,\ldots,p_n) \in Supp(H^i(\tilde{R})) \mid {\rm\
>some\ } p_j < - \ell \}.$$
\Box

%Date: Fri, 6 Mar 1998 14:09:01 -0500
%From: Mike Stillman - Math Prof <mike@math.cornell.edu>

I talked to Lou Billera briefly about the location of
sign placements.  He mentioned the following trivial
argument.

Let A be the 2 by n matrix of the locations of the edges in the fan
(in the plane), in order.  The matrix [c d]A gives the linear combinations
that we were looking at.  The question is, what are the possible
sign arrangements of the resulting vector?  The vector [c d]
determines a line perpendicular to [c d].  A "+" in a
component means that the corresponding edge
is on the positive side of this line, a "-" means that it is on the
other side.  So the possible sign arrangements are obtained by seeing which
edges are on each side of this line.  Clearly all the "-"'s must be
contiguous, from this picture.  This also handles the case of "0"'s,
since the 0's must be on this line.

%%%%%%%Mike's original (wrong??) proof from section 2
\corollary{coarse image} If $I \in \Sigma_i$, then the projection
$\phi(C_I)$ is a pointed cone; that is, if $x,-x\in\phi(C_I)$ then
$x=0$.

\proof If $\alpha,\beta$ are preimages of $x,-x$ then
$\alpha+\beta\in(\ker \phi)\cap C_I$.
But
$$
H^i(\O_X) = \bigoplus_{\alpha \in \ker\phi} H^i_*(\O_X)_\alpha = 0,
$$
unless $i=0$, in which case it is supported in degree $\alpha=0$.
Thus
if $\alpha \in \ker \phi$, then $\alpha$ is not in the support of
$H^i_*(\O_X)$. Therefore if $\alpha \in C_I \cap \ker \phi$, then $\alpha = 0$.
\Box

####################################################
Some stuff I might want to keep... (MES)

In this section we describe algorithms for computing cohomology of sheaves on
toric varieties.  We regard as input to these algorithms the data which give
the toric variety:

\item{(a)} The maps
  $$\ZZ^d \longrightarrow \ZZ^n \longrightarrow Ch(X) \longrightarrow 0,$$

\item{(b)} The homogeneous coordinate ring $R = k[x_1, \ldots, x_n]$, where
$\deg(x_i) = \pi(e_i)$, with $e_i$ the $i$th unit vector.

\item{(c)} The irrelevant ideal $B$ in $R$.

\bigskip

We regard as known algorithms for computing $Ext^i(B^{[\ell]},M)$, the degree
$e \in Ch(X)$ part of a graded module, and also
the bounding hyperplanes for a convex cone.

The main algorithm presented here is to compute functions $f_i$ from
finitely generated graded $R$-modules to $\ZZ_{\geq 0}$, which satisfies the
following property:
  For all $\ell \geq f_i(M)$, the degree zero part of the map
  $$Ext^i(B^{[\ell]},M) \longrightarrow H^i_*(\tilde(R))$$
is an isomorphism.
Again, we may look at this from the opposite point of view.

\theorem{}
Let $P$ be a finitely generated graded $S$ module ( with the coarse
grading). Let $F_{\bullet}$ be a minimal free resolution of $P$ with
$$F_i=\bigoplus_{\alpha\in\Ch(X)}S(-\alpha)^{\beta_{i,\alpha}}.$$
Then the set of $D \in Ch(X)$ such that
$$\Ext^i_S(B^{[\ell]}, P)_D\longrightarrow H^i_*(\tilde{P})_D$$
is an isomorphism contains the complement of the union of the
translated pointed cones
  $$\alpha + \pi(p_I - \ell e_j) + \pi(C_I),$$
for all $q \geq 0$, all $I \in \Sigma_q$, all $j \in I$, and all
$\alpha \in Ch(X)$ such that $\beta_{q-i+1,\alpha}\neq 0$ or
$\beta_{q-i,\alpha}\neq 0$.
\bigskip